\documentclass{amsart}

\usepackage[protrusion=true
]{microtype}
\usepackage[foot]{amsaddr}
\usepackage{amssymb}
\usepackage{booktabs}
\usepackage{tikz}\usetikzlibrary{positioning}
\usepackage{natbib}
\renewcommand{\cite}{\citep}
\newtheorem{scheme}{Argumentation Scheme}

\title{Evidence, Proofs, and Derivations}
\author[Andrew Aberdein]{Andrew Aberdein$^{*}$}
\address{$^{*}$School of Arts \& Communication, Florida Institute of Technology, Melbourne FL.}
\date{\today}
\thanks{Forthcoming in \textit{ZDM mathematics education} {51(4)}, 2019.}
\begin{document}
\begin{abstract}
The traditional view of evidence in mathematics is that evidence is just proof and proof is just derivation. There are good reasons for thinking that this view should be rejected: it misrepresents both historical and current mathematical practice. Nonetheless, evidence, proof, and derivation are closely intertwined. This paper seeks to tease these concepts apart. It emphasizes the role of argumentation as a context shared by evidence, proofs, and derivations. The utility of argumentation theory, in general, and argumentation schemes, in particular, as a methodology for the study of mathematical practice is thereby demonstrated. Argumentation schemes represent an almost untapped resource for mathematics education. Notably, they provide a consistent treatment of rigorous and non-rigorous argumentation, thereby working to exhibit the continuity of reasoning in mathematics with reasoning in other areas. Moreover, since argumentation schemes are a comparatively mature methodology, there is a substantial body of existing work to draw upon, including some increasingly sophisticated software tools. Such tools have significant potential for the analysis and evaluation of mathematical argumentation. The first four sections of the paper address the relationships of evidence to proof, proof to derivation, argument to proof, and argument to evidence, respectively. The final section directly addresses some of the educational implications of an argumentation scheme account of mathematical reasoning.
\end{abstract}
\maketitle

\section{Evidence versus Proof}\label{EvsP}
The traditional view of evidence in mathematics is that evidence is just proof. 
Donald Martin poses the question,
`What then does count as mathematical evidence? There is, of course, an obvious answer to the question of how one can come to know the truth of a mathematical proposition: namely, proof. Indeed, this may seem the \emph{only} way to establish mathematical truth' \cite[216]{Martin98}.
However, as Martin acknowledges, there are good reasons for not accepting this view. Crucially, it does not comport with actual mathematical practice. 
Here, for example, is a short, melancholy narrative of an unsuccessful proof attempt, unusual only in being so self-contained. The mathematician Va\v{s}ek Chv\'atal proposed as a conjecture a generalization of the well-known Sylvester--Gallai theorem. Unable to supply a proof, he comments that `we present meagre evidence in support of this rash conjecture' \cite[175]{Chvatal04}.
However, the last lines of his published paper read:
`\emph{Received July 14, 2002, and in revised form April 14, 2003. Online publication December 31, 2003. Note added in proof}. In September 2003 Xiaomin Chen proved Conjecture 3.2' \cite[195]{Chvatal04}.
Conjecture 3.2 is Chv\'atal's `rash conjecture'; where Chv\'atal had evidence, Chen had proof. If proof were the only evidence there is in mathematics, then either Chv\'atal would have had nothing to say, or Chen nothing to add.

Here is a more protracted exhibition of the use of evidence beyond proof in mathematical practice.
One of the last publications of the celebrated mathematician (and Nobel laureate) John Nash was a collection of essays on open problems, edited in collaboration with Michael Rassias.
The editors remark of their choice of problems that
`Some were chosen for their undoubtable importance and applicability, others because they constitute intriguing curiosities which remain unexplained mysteries on the basis of current knowledge and techniques, and some for more emotional reasons' \cite[vi]{Nash16}.
Their selection includes many of the best known open problems in mathematics, each discussed by a leading expert. Not every essay explicitly references the evidence for these conjectures, but many do. Here is a sample (boldface emphasis mine throughout; internal citations omitted):
\begin{description}
\item[$\mathsf{P}$ versus $\mathsf{NP}$]
`Of course, this gives even more dramatic \textbf{evidence} that \textsc{GraphIso} is not \textsf{NP}-complete: if it was, then \emph{all} \textsf{NP} problems would be solvable in $n^{\mathrm{polylog}\,n}$ time as well' \cite[20]{Aaronson16}.
\item[Montgomery's Pair Correlation Conjecture]
`The agreement with the first million zeros is poor, but the agreement near zero number $10^{12}$ is close, near perfect near zero number $10^{16}$, and even better near zero number $10^{20}$. These results provide massive \textbf{evidence} for Montgomery's conjecture' \cite[155]{Barrett16}.
\item[The Birch--Swinnerton-Dyer Conjecture]
`Today, the numerical \textbf{evidence} in support of both the weak and full Birch--Swinnerton-Dyer conjecture is overwhelming, and probably more extensive than for any other conjecture in the history of mathematics' \cite[213]{Coates16}.
\item[The {{E}rd\H{o}s}--{S}zekeres Problem]
`The special cases $A=I$ and $A=\{i_{1},i_{n}\}$ correspond to minimal \textbf{evidence} for existence of an $n$-cup and an $n$-cap' \cite[363]{Morris16}.
\item[The Hadwiger--Nelson Problem]
`Ronald L. Graham \dots\ 
cites a theorem of Paul O'Donnell (see [34, 49]) showing the existence of 4-chromatic unit-distance graphs of arbitrarily large girth (Theorem 28 below) as ``perhaps, the \textbf{evidence} that $\chi$ is at least 5.''' \cite[441]{Soifer16}.
\item[Erd\H{o}s's Unit Distance Conjecture]
`There is no strong \textbf{evidence} supporting the assumption that such an example cannot exist' \cite[468]{Szemeredi16}.
\item[Goldbach's Conjecture]
This is the strongest \textbf{evidence} we have that for even $n$, $r_{2}(n) \sim n\mathfrak{S}_{2}(n)$' \cite[486]{Vaughan16}.
\item[The Hodge Conjecture]
`The following theorem proved in [7] is the best known \textbf{evidence} for the Hodge conjecture
' \cite[538]{Voisin16}.
\end{description}
None of these examples of evidence refers to proof; indeed, in many of them the evidence is for a conjecture notable for being unproven. Nonetheless, we may see that mathematical evidence can be dramatic, massive, overwhelming, strong---or, conversely, minimal. But even overwhelming evidence falls short of proof.
Notice also that evidence can diverge from proof on more than one dimension. Some of these cases fall short of proof since they are proofs of some weaker conjecture that provides bounds on the headline conjecture, or is otherwise suggestive of its truth (e.g. Hadwiger--Nelson, Hodge); others differ in kind, rather than degree, by offering numerical, that is empirical, evidence  (e.g. Montgomery, Birch--Swinnerton-Dyer). Nor are considerations of evidence short of proof confined to the contemplation of stubbornly unsolved problems. On the contrary, mathematicians must often consider such evidence when making decisions about their own careers and those of their students. Crucially, as James Franklin notes, every Ph.D. supervisor must ensure that the problems their students tackle are open, but likely to yield to a few years work \cite[2]{Franklin87}. Similar considerations of evidence bear on grant awards, tenure decisions, and many other administrative aspects of mathematical practice.

\section{Proof versus  Derivation}\label{PvsD}
The `standard view' of the relationship between informal, or everyday proof and derivation, its formal counterpart, is that
`informal proofs are just sloppy, incomplete versions of formal proofs' \cite{Pawlowski17}.
On this account, the only proofs worthy of the name would either already be derivations or be rewritable as derivations in a more or less trivial manner.
Nonetheless, it has lately become 
`a common observation that the proofs that mathematicians write on blackboards and publish in journals are not like the derivations that appear as objects in proof theory' \cite[401]{Larvor16}.
However, it can be a challenge to pin down exactly where the divergence arises.

Before going any further, I should clarify what I am \emph{not} saying.
Nothing in this paper is intended to contradict what is sometimes termed the \emph{Formalizability Thesis}: that every proof can (in principle, at least) be formalized as a derivation.%
\footnote{This is sometimes referred to as \emph{Hilbert's Thesis}, although that name is more properly reserved for the narrower claim that every proof can be formalized as a derivation in first-order logic \cite{Kahle18}.}
Crucially, the Formalizability Thesis is an existence claim: it states that for every proof there is a derivation. It makes no claim as to whether the derivation is known, or accessible, or surveyable.
It does not require that the internal structure of a proof bear any resemblance to that of the corresponding derivation.%
\footnote{This reflects what has been called \emph{Tait's Maxim}: `The notion of formal proof was invented to study the existence of proofs, not methods of proof' \cite[114]{Baldwin13}.}
It also makes no claim as to the metaphysical status of the relationship between the proof and the derivation.
For present purposes it suffices to observe that accepting that there is a derivation for every proof is consistent with proofs and derivations being profoundly different entities and specifically does not entail that knowing a proof has any connection to knowing a derivation.

The proof/derivation distinction has been drawn in multiple, conceptually distinct ways. I shall survey some of the most important---semantic/syntactic; normative/theoretical; act/object---before addressing how they are related. 
One influential account of the proof/derivation distinction situates it in the different manner of expression of proofs and derivations. Thus Yehuda Rav, in a highly influential paper, declares that he will
`understand by \emph{proof} a conceptual proof of customary mathematical discourse, having an irreducible semantic content, and distinguish it from \emph{derivation}, which is a syntactic object of some formal system' \cite[11]{Rav99}.
Likewise, Jody Azzouni draws the distinction between
`formal derivations, which occur in artificial languages, and mathematical proof, which occurs in natural languages' \cite[247]{Azzouni13}.
Analogously, Keith Weber and Lara Alcock
`define a syntactic proof production to occur when the prover draws inferences by manipulating symbolic formulae in a logically permissible way' and `define a semantic proof production to occur when the prover uses instantiations of mathematical concepts to guide the formal inferences that he or she draws' \cite[209]{Weber04}.
While all of these points are correct, we may feel that they do not strike to the heart of the distinction.

By contrast, Gila Hanna draws the distinction in a manner that brings out the social aspects of the contrast:
\begin{quotation}
\begin{enumerate}
\item Formal proof: proof as a theoretical concept in formal logic (or metalogic), which may be thought of as the ideal which actual mathematical practice only approximates. 
\item Acceptable proof: proof as a normative concept that defines what is acceptable to qualified mathematicians \cite[6]{Hanna90}.
\end{enumerate}
\end{quotation}
The mathematician John Baldwin, asserts that he takes Hanna's `formal/acceptable to be the same distinction' \cite[71]{Baldwin16} as that he draws between `two degrees of formalization', respectively into a formal syntax and rules of inference or within (sufficiently mathematical) natural language \cite[89]{Baldwin13}.
So, for Baldwin at least, Hanna's distinction is coextensive with that between proofs and derivations. 
However, Hanna stresses the essentially social role of acceptability in proof which has no counterpart in derivation.
Authors who emphasize the argumentational nature of proof, to which we will return below, specifically draw attention to this factor. Thus 
Todd CadwalladerOlsker states that the
`formal view of proof is contrasted with the view of proofs as arguments intended to convince a reader' \cite[33]{CadwalladerOlsker11}
and Trevor Bench-Capon stresses that
`argumentation is \dots\ an activity which has to be actively engaged with, whereas a proof is an object to be understood and admired' \cite{BenchCapon12}.
The latter observation points towards an aspect of the proof/derivation contrast that Hanna's distinction overlooks: derivations are objects; proofs are acts.

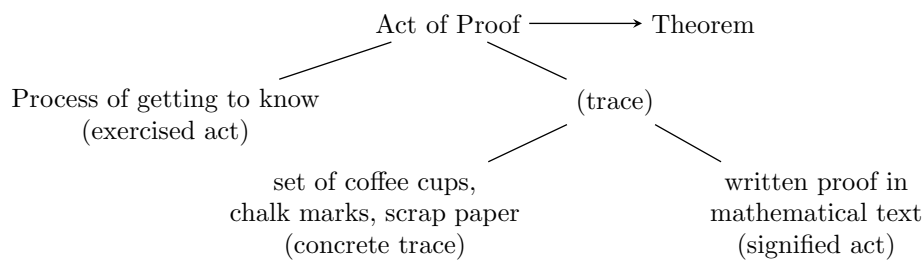
\begin{figure}[htbp]
\begin{center}
\begin{tikzpicture}[>=stealth,line width=.5pt,every node/.style={align=flush center}] 
\node (Act) {Act of Proof};
\node (Thm) [right=1.5 of Act] {Theorem};
\node (Proc) [below left=.7 of Act] {Process of getting to know\\ (exercised act)};
\node (trace) [below right=.7 of Act] {(trace)};
\node (conc) [below left=.7 of trace] {set of coffee cups,\\ chalk marks, scrap paper\\(concrete trace)};
\node (sign) [below right=.7 of trace] {written proof in\\ mathematical text\\ (signified act)};
\draw[->]	(Act) -- (Thm);
\draw[-]	(Proc) -- (Act);
\draw[-]	(trace) -- (Act);
\draw[-]	(trace) -- (conc);
\draw[-]	(trace) -- (sign);
\end{tikzpicture}
\caption{G\"oran Sundholm's act/process/trace distinction \cite[after][948]{Sundholm12a}}
\label{Goran}
\end{center}
\end{figure}

In support of an account of proofs as acts, Joseph Goguen argues that the only sort of proofs `that can actually happen in the real world are \emph{proof events}, or \emph{provings}, which are actual experiences, each occurring at a particular time and place, and involving particular people, who have particular skills as members of an appropriate mathematical community' (\citealp{Goguen01}; \citealp[see also][489 ff.]{Stefaneas12}).
Goguen is not alone in drawing the distinction this way. Oswaldo Chateaubriand also distinguishes between `provings and idealized proofs' \cite[41]{Chateaubriand03}, while tracing the temporal conception back to Brouwer.
Goguen also observes that provings often have internal temporal structure---their components must be executed in the right order---and are thus ultimately {proof processes}:
\begin{quotation}
The efficacy of some proof events depends on the components of a proof object being seen to be given in a certain temporal order, e.g., Euclidean geometric proofs, and commutative diagrams in algebra; in some cases, the order may not be easily infered from just the diagram. Therefore we must generalize from static proof objects to \emph{proof processes}, such as diagrams being drawn, movies being shown, and Java applets being executed \cite{Goguen01}.
\end{quotation}
Lest this conception of proof seem to have drifted too far from what is conventionally labelled as proof, it is important to reflect on what is meant by an act. 
G\"oran Sundholm draws a useful distinction between an act, the subjective process that comprises the act, and the trace that the act leaves behind \cite[948]{Sundholm12a}. In the case of proof, we may further distinguish the concrete traces (marks on the blackboard, empty coffee cups, and so on) from the informational trace, a written proof (see Fig.~\ref{Goran}). The latter signifies the act and serves as a blueprint or recipe whereby the act may be repeated. Just as recipes aren't very nourishing unless you carry them out, likewise the written proof only works as a proof if it is actually carried out.
Although the contrast is more stark in the case of proofs, we may draw a similar contrast between derivations and derivation traces: the derivation is a mathematical object; its trace is a written counterpart to that object (a formal proof, or a sequence of code, or the like).

In summary, we shall say that proofs make irreducible use of the semantics of natural language; involve a normative appeal to acceptance by a mathematical audience; and are events that involve participants and extend over time. 
Conversely, derivations are mathematical objects that represent formal relationships between propositions in an artificial language.
Of course, this definition of proof is stipulative. If we treat the three contrasts discussed above---semantic/syntactic; normative/theoretical; act/object---as separate axes, then my definitions of proof and derivation occupy diagonally opposite vertices: $\{0,0,0\}$ and $\{1,1,1\}$, as it were. That leaves six other, unlabelled vertices. I shall not attempt a complete taxonomy here, but some conceptions of proof rival to that defended here may be seen to correspond to some of these other vertices.

\section{Argument versus Proof}\label{AvsD}
We saw in \S\ref{EvsP} that mathematicians make use of a concept of evidence distinct from proof and we saw in \S\ref{PvsD} that the practice of mathematical proof cannot be reduced to derivation. But nor can we discard derivations: they are mathematical objects in good standing (whatever that means---I make no ontological claim), and acknowledged as such by the mathematicians who study them. 
Furthermore, at least in the eyes of many mathematicians, it is the existence of a derivation that is the ultimate guarantor of the truth of a theorem.\footnote{For example, for Saunders Mac Lane, ‘the test for the correctness of a proposed proof is by formal criteria and not by reference to the subject matter at issue’ \cite[378]{MacLane86} and Thomas Hales characterizes formal proof as providing ‘a thorough verification of my own research that goes beyond what the traditional peer review process has been able to provide’ \cite[1378]{Hales08}.}
So, if we are to respect actual mathematical practice, we must accommodate all three concepts: evidence, proof, and derivation. Is there a common framework into which they may be subsumed?
One candidate is argument. I have already alluded to an argumentation-based account of mathematical reasoning. In its barest form, this states that
a mathematical argument is `a more liberalized version of the notion of mathematical proof' \cite[157]{VanBendegem05}. That is, proofs are a special case of mathematical argument.
If these arguments are understood as similar in kind to the arguments of non-mathematical discourse, it follows that an account of mathematical reasoning may be given using the tools of argumentation theory.

Elsewhere I have referred to proofs$^{*}$, where the asterisk indicates that the proof does not command universal assent \cite[2]{Aberdein09b}.
These include unsurveyably long proofs$^{*}$, diagrammatic proofs$^{*}$, proofs$^{*}$ that depend on contested axioms, computer-assisted proof$^{*}$, experimental proof$^{*}$, proof$^{*}$ by probabilistic methods, and so on.
While it is controversial in each case whether proofs$^{*}$ should lose the asterisk
\citetext{for relevant discussion, see \citealp{Fallis97,VanBendegem05,Paseau15}},
each of these is uncontroversially a mathematical argument.
Some proofs$^{*}$ also uncontroversially correspond to derivations. 
For others this is less obvious; their admission as proofs would pose a challenge to the Formalizability Thesis discussed in \S\ref{PvsD}. 
The Formalizability Thesis, that `every proof can be formalized as a derivation', is one clause of what has been called \emph{Leibniz's Thesis}; the other clause states that `every acceptable argument of (informal) mathematics is a proof' \cite[17]{Berk82}.
This implies that mathematical arguments that are not proofs are unacceptable. But what does it mean for a mathematical argument to be acceptable? Interpreted widely, as `acceptable in some mathematical practice', then the thesis is demonstrably false; interpreted narrowly, as `acceptable as proof', then the thesis seems almost tautological.
Nonetheless, as we have acknowledged above, being acceptable to a mathematical audience is a necessary requirement for a proof. Below I shall sketch what that might mean.

Arguments resemble proofs, at least as defined above, much more closely than they do derivations. Both arguments and proofs are acts (or processes comprised of acts) whereas derivations are objects. Yet, as we saw in the last section, proofs and derivations both leave traces: written blueprints or recipes. Likewise, we can consider argument traces.
The passing resemblance between a proof trace and a derivation trace explains how two such dissimilar things as proofs and derivations ever came to be conflated.
A derivation is a directed graph. Its nodes are (formal counterparts of) truths of mathematics and its edges represent logical deductions. A determination of whether a derivation is sound requires a choice of axioms and of logical system. Relative to that choice, the derivation is sound if all its source nodes (those which do not depend on other nodes) are axioms and all its edges are valid.
Proof traces (and argument traces in general) are also directed graphs. The nodes of a proof trace are mathematical truths, expressed in a suitably augmented natural language, and the edges are arguments of some kind.

So, the same theorem will be linked to (at least) two structures: a proof trace and a derivation, each of which is comprised of directed graphs. How do these two structures relate to each other?
In principle, any node in either structure could be a node in the other, at least assuming the Formalizability Thesis holds. In practice, there is rather less overlap. Derivations are substantially more verbose than proofs: they contain many more intermediate statements and they take everything back to axioms, as proofs characteristically do not. What's more, they may follow an entirely different path from any recognisable proof. (Recall Tait's Maxim, cited above.) And, for that reason, proof traces may contain nodes that do not figure in any corresponding derivation. But at least the nodes of the two structures are held to the same standard: in each case, they must be true.
The edges, on the other hand, must be judged differently. As deductions in a system of formal logic, the edges of a derivation can be assessed by the canons of that system: do they instantiate an admissible rule of inference? Of course, the edges of the proof structure may also be logical deductions, in which case they can be held to that standard too. But, with some rarefied exceptions, mathematical proofs are seldom purely logical proofs. So the edges of the proof structure must be judged instead by standards appropriate to arguments of the relevant kind; that is, by the standards of (some suitably localized form of) argumentation theory.

Multiple methodologies for the analysis and evaluation of informal arguments have been proposed, but here I shall focus on the method of argumentation schemes. This is an ancient idea in origin, deriving from the topoi or loci of classical rhetoric, but it has been reinvigorated in recent years. 
An argumentation scheme is a template that captures a stereotypical pattern of reasoning. 
Different schemes are fine-tuned to capture the idiosyncrasies of different types of argument---and can be brought to bear to determine when they have been used cogently.
This evaluative function is largely the role of the critical questions, which most schemes contain. If in some instantiation of a scheme these questions cannot be adequately addressed, then the argument fails. Deductive inference rules can be understood as argumentation schemes, but the method comes into its own when applied to non-deductive reasoning. 
Douglas Walton and colleagues characterize many such schemes as special cases of the very general scheme that they call Defeasible Modus Ponens \cite[366]{Walton08}.
In Scheme~\ref{DMP}, I have adapted their presentation of this scheme to bring out the resemblance between argumentation schemes and Toulmin layouts \cite[101]{Toulmin58}:%
\footnote{For a more protracted discussion of how these two models of reasoning are related, see \cite[28 ff.]{Pease11}. For an alternative account, see \cite[1070]{Konstantinidou13}.}

\begin{scheme}\label{DMP}
Defeasible Modus Ponens
\end{scheme}
\begin{description}
\item[Data]$P$.
\item[Warrant] As a rule, if $P$, then $Q$.\\
Therefore, \dots 
\item[Qualifier] presumably, \dots
\item[Conclusion] \dots\ $Q$. 
\end{description}
\begin{center}
Critical Questions
\end{center}
\begin{enumerate}
\item \textbf{Backing:} What reason is there to accept that, as a rule, if $P$, then $Q$?
\item \textbf{Rebuttal:} Is the present case an exception to the rule that if $P$, then $Q$?
\end{enumerate}
In principle, any argument that can be represented using Toulmin layouts could be represented using Scheme~\ref{DMP}. However, the great strength of the argumentation scheme methodology lies in its diversity: there are many more schemes to choose from. One of the most extensive surveys of general purpose schemes distinguishes more than 90 different types \cite[308 ff.]{Walton08}. Some of these off-the-shelf schemes are directly applicable to mathematics, but yet more can be produced to order.

Elsewhere I draw a distinction between three classes of scheme that are of immediate relevance to proof \cite[366 f.]{Aberdein13a}:
\begin{itemize}
\item {\bf A-schemes} {correspond directly to derivation rules}. 
(Equivalently, we could think in terms of a single A-scheme, the `pointing scheme' which picks out a  derivation whose premisses and conclusion are formal counterparts of its data and claim.) 
\item {\bf B-schemes} are {exclusively mathematical arguments}: high-level algorithms or macros. Their instantiations {correspond to substructures of derivations} rather than individual derivations (and they may appeal to additional formally verified propositions). 
\item {\bf C-schemes} are even looser in their relationship to derivations, since {the link between their data and claim need not be deductive}. 
Specific instantiations may still correspond to derivations, but there will be no guarantee that this is so and no procedure that will always yield the required structure even when it exists.
Thus, where the qualifier of A- and B-schemes will always indicate deductive certainty, the qualifiers of C-schemes may exhibit more diversity. Indeed, different instantiations of the same scheme may have different qualifiers.
\end{itemize}
We are now in a position to analyze what it means for a proof to be acceptable. For a given mathematical audience, that is a community of mathematicians who share the same standards, a proof will count as acceptable if the schemes it instantiates are ones which that audience judges consistent with mathematical rigour. 
In other words, the audience must be convinced that each step of the proof instantiates a scheme in such a way that it earns the qualifier `rigorously' (or stronger).
Exactly which schemes are deemed consistent with rigour will vary by audience. Different areas of mathematics can have somewhat different standards and even within the same area there are different audiences: the audience for a research article is not identical to that for an undergraduate lecture.

Nonetheless, we may generalize. For most audiences, A-schemes are invariably admissible. Likewise, most audiences will admit many B-schemes. However, some B-schemes can be highly complex. As such they may only be admissible to a narrow audience: professionals in a particular subspecialty, say. Such schemes may need to be broken into simpler steps (probably also instantiating simpler B-schemes) for consumption by a wider audience of research mathematicians, let alone for a student audience.
Indeed, in some pedagogic contexts, it may be expedient to substitute a less (than) rigorous C-scheme for an intricate B-scheme, thereby ‘handwaving’ through an aspect of the proof unsuitable for a given audience. C-schemes themselves pose a further challenge. Since not all their instantiations are deductive, it may be tempting to assume that they have no place in rigorous argument, at least as most mathematical audiences understand it. Even were this so, they would still be welcomed by some mathematical audiences, notably those that see no need for the asterisks on some of the more outr\'e species of proof$^{*}$. However, in practice, even quite conservative mathematical audiences can find some C-schemes admissible. 

\begin{table}[htbp]
\begin{center}
\caption{Summary of reasoning types and schemes \cite[247]{Aberdein13}}\label{Schemes}
\begin{tabular}{l}
\toprule
\begin{minipage}{.85\linewidth}
\begin{enumerate}
\item Reasoning
\begin{enumerate}
\item Retroduction
\begin{enumerate}
\item Argument from Gradualism
\item Argument from Positive Consequences
\item Argument from Evidence to a Hypothesis
\end{enumerate}
\end{enumerate}
\item Source Based
\begin{enumerate}
\item Citation
\begin{enumerate}
\item Appeal to Expert Opinion
\item Argument from Danger 
\end{enumerate}
\item Intuition
\begin{enumerate}
\item Argument from Position to Know
\end{enumerate}
\item Meta-Argument
\begin{enumerate}
\item Ethotic Argument
\end{enumerate}
\item Closure
\begin{enumerate}
\item Argument from Ignorance
\end{enumerate}
\end{enumerate}
\item Rule Based
\begin{enumerate}
\item Generalization
\begin{enumerate}
\item Argument from Example
\end{enumerate}
\item Definition
\begin{enumerate}
\item Argument from Definition to Verbal Classification
\end{enumerate}
\end{enumerate}
\end{enumerate}
\end{minipage}\\
\bottomrule
\end{tabular}
\end{center}
\end{table}

Some of the general purpose C-schemes that may be accepted by some mathematical audiences are listed in Table~\ref{Schemes} 
This is by no means an exhaustive list; nor is the classification of schemes it presumes unarguable. In the article from which Table~\ref{Schemes} is taken, I discuss the application to informal mathematical reasoning of ten schemes drawn from \cite{Walton08}. I only have space to address a couple of them in comparable detail here, but I will briefly survey the rest. I classified three schemes as retroduction, that is reasoning backwards from a conclusion, a common strategy in mathematical argument since antiquity. Argument from Gradualism comprises multi-step argumentation. This can be fallacious when the steps are individually dubious or collectively improbable, yielding a slippery slope fallacy. Conversely, if all the steps are deductively valid, it would correspond to a derivation. Some instances of this scheme, however, are neither fallacious nor deductively valid: their admissibility will turn on the standards of rigour of the audience. Argument from Positive Consequences is well described by its name. It finds a place in mathematics in the informal justification of axioms and hypotheses.\footnote{The mathematical uses of Argument from Positive Consequences are also discussed, together with some other schemes not in Table 1, in work by Nikolaos Metaxas and colleagues \citetext{\citealp[84]{Metaxas15}; \citealp[387]{Metaxas16}}.} I will discuss Argument from Evidence to a Hypothesis at greater length in the next section.

Several of the schemes in Table~\ref{Schemes} are source-based: that is, they involve reasoning to a mathematical conclusion from a source external to the reasoning process. The most familiar of these reasoning patterns may be the use of citation, which is ubiquitous in most mathematical research. Appeal to Expert Opinion is an obvious fit for such reasoning, and also one of the most widely studied schemes in the argumentation scheme literature. Argument from Danger, which I also included under citation, is of much narrower application. It captures arguments to the effect that one should not act in a way that would be a danger (to oneself or others). I argued that this is employed in what is sometimes (facetiously) called ‘proof by intimidation’: dismissing objections as trivial, for example. Argument from Position to Know describes a pattern of argument where an informant is taken to have a privileged source of information. The obvious non-mathematical example would be an eye witness; in a mathematical context, this scheme might be employed to characterize appeals to intuition, which some authors analogize to perception \cite[for example,][]{Chudnoff13c}. Ethotic Arguments are appeals to the character (or ethos) of some individual. It has been shown empirically that such factors can also influence how mathematical results are received \cite{Inglis09a}. Argument from Ignorance is generally treated as a fallacy, since not knowing that something is false is not usually a legitimate reason for treating it as true. Yet it can be a reliable inference in special circumstances, some of which can be found in mathematical reasoning, such as reporting the outcome of an exhaustive search.

The last group of schemes in Table~\ref{Schemes} are those involving the justification or application of rules. Argument from Definition to Verbal Classification could potentially correspond to a step in a derivation, but it can also be used more casually, and is vulnerable to outright misuse, for example when a definition is subtly misstated. For a more considered example, consider the following:
\begin{scheme}
Argument from Example\label{Example}
\end{scheme}
\begin{description}
\item[Premise] In this particular case, the individual $a$ has property $F$ and also property $G$.
\item[Conclusion] Therefore, generally, if $x$ has property $F$, then it also has property $G$.
\end{description}
\begin{center}
Critical Questions
\end{center}
\begin{enumerate}
\item Is the proposition claimed in the premise in fact true?
\item Does the example cited support the generalization it is supposed to be an instance of?
\item Is the example typical of the kinds of cases the generalization covers?
\item How strong is the generalization?
\item Do special circumstances of the example impair its generalizability? \citep[314]{Walton08}
\end{enumerate}
Many instantiations of Argument from Example will not be acceptable at least to modern audiences. Enumerative induction, for instance, can be characterized in terms of this scheme. But there are other cases where Argument from Example meets modern standards of rigour. In order to do so, adequate answers to all the critical questions will need to be provided. In particular, question (4) will require an exceptional answer: one justifying the revision of the qualifier from `generally' to `necessarily'.%
\footnote{For a more extensive discussion of mathematical uses of Scheme~\ref{Example}, see \cite[244]{Aberdein13}.}

\section{Argument versus Evidence}\label{AvsE}
We have explored the relationship between proofs and derivations and between each of these and arguments. We have not yet addressed the relationship between evidence and argument, but alert readers may have spotted the following scheme in Table~\ref{Schemes}:
\begin{scheme}\label{Hypothesis}
Argument from Evidence to a Hypothesis
\end{scheme}
\begin{description}
\item[Major Premise] If $A$ (a hypothesis) is true, then $B$ (a proposition reporting an event) will be observed to be true.
\item[Minor Premise] $B$ has been observed to be true, in a given instance.
\item[Conclusion] Therefore, [presumably,] $A$ is true.
\end{description}
\begin{center}
Critical Questions
\end{center}
\begin{enumerate}
\item Is it the case that if $A$ is true, then $B$ is true?
\item Has $B$ been observed to be true?
\item Could there be some reason why $B$ is true, other than its being because of $A$ being true? \citep[331 f.]{Walton08}
\end{enumerate}
This is the pattern of reasoning called abduction by Charles Peirce (for a direct comparison, see \citealp[22]{Pease11}; for an alternative approach to abduction in terms of argumentation schemes, see \citealp[386]{Metaxas16}). 
To observe how this scheme might work in practice, consider the following example of evidence for Goldbach's Conjecture (GC), as reconstructed by the philosopher Alan Baker:
\begin{quotation}
Another line of auxiliary argument might be based on the various partial results relating to GC. In 1931, Schnirelmann proved that every even number can be written as the sum of not more than 300,000 primes(!).
This upper bound on the number of primes required has since lowered to 6 (Ramar\'e 1995). In addition, Chen (1978) proved that all sufficiently large even numbers are the sum of a prime and the product of two primes. Such results do not seem to make the truth of GC any more likely. But perhaps they provide evidence that GC is provable \cite[71]{Baker07}.
\end{quotation}
Here the hypothesis $A$ is the proposition that GC is provable and the `event' (or events) that $B$ reports on are a series of proofs of weaker, but related, conjectures. As Baker acknowledges, this is at best suggestive. Critical question (3) is completely open and the argument would only sustain a very weak qualifier.
On the other hand, some of the examples of argument from evidence cited in \S\ref{EvsP} are much more convincing.
For instance, the argument for Montgomery's Pair Correlation Conjecture can easily be reconstructed as an instance of Scheme~\ref{Hypothesis}, with $A$ a statement of the conjecture and $B$ the observation that `the agreement near zero number $10^{12}$ is close, near perfect near zero number $10^{16}$, and even better near zero number $10^{20}$' \cite[155]{Barrett16}.
Of course, as stressed in \S\ref{EvsP}, the authors of this argument are not claiming that this is a proof. Hence this instance of Scheme~\ref{Hypothesis} is not one that they would treat as acceptable.

Are there any instances of Scheme~\ref{Hypothesis} which a mainstream mathematical audience would treat as acceptable? One possible affirmative answer brings us back to derivations. 
In his discussion of the relationship between proofs and derivations, Richard Epstein makes the following bold assertion:
`A proof in a fully formal system of logic that a claim follows from some axioms is not a proof in mathematics. It is evidence that can be used in a mathematical proof' \cite[274]{Epstein12}.
In the terminology of this paper, Epstein may be read as saying that a derivation, or more properly, a derivation trace, should not be mistaken for a proof (or proof trace), but can be employed in a proof as evidence. On the account of proof defended above, such employment would require use of an argumentation scheme, presumably along the lines of Scheme~\ref{Hypothesis}.
For proofs where the corresponding derivation is comparatively short and straightforward, this may seem to be an unnecessarily scrupulous point. Yes, in principle, the derivation is a mathematical object, but its trace will be similar enough to a proof trace for it to be treated as such with comparatively little additional effort, at least for audiences familiar with such things. In these sorts of cases, the more characteristic response to a derivation trace would be to `reverse engineer' it into a proof in this manner.
But, for longer or more technical derivations, this process is less practical. In such cases, Epstein's analysis seems correct. If there is a sufficiently strong reason for accepting the derivation trace as reliable, then this would be an acceptable instance of Scheme~\ref{Hypothesis}.

Unsurveyable computer-assisted proof$^{*}$ represents a similar application of Scheme~\ref{Hypothesis}. Proofs such as Kenneth Appel and Wolfgang Haken's proof of the Four Colour Conjecture or Thomas Hales's proof of Kepler's Conjecture depend essentially on unsurveyably vast computer calculations. Within the human readable component of the proof, these calculations play an evidential role. Thus, if these proofs are acceptable to the mathematical community, then the steps introducing the computer calculations can be seen as acceptable instances of Scheme~\ref{Hypothesis}.
Of course, while the consensus seems now to be in their favour, I have flagged such results as proof$^{*}$ because their acceptability remains a topic of debate.
However, if the existence of a derivation is what ultimately underwrites the acceptance of a theorem, then a mechanically verified derivation trace may be an (even more) acceptable instance of Scheme~\ref{Hypothesis}. Successful formalization projects, such as Georges Gonthier's work on the Four Colour Theorem or Hales's Flyspeck Project, have provided just such corroboration \cite{Gonthier08a,Hales17}.%
\footnote{It may be objected that this results in a regress, since the software checking the derivation  trace must itself be checked. However, it is what Hales has called `a rather manageable regress' \cite[1376]{Hales08}. The kernel of such proof checking software is very carefully designed to be small enough and clear enough to be amenable to thorough human checking.}

\section{Conclusions for Education}\label{Education}
One of the morals of this paper is that it is important to carefully distinguish proof from derivation. This observation certainly has profound implications for the teaching of mathematics. 
However, it is scarcely novel: I have already cited a paper from almost thirty years ago which addresses the educational implications of the distinction \cite{Hanna90}.
The role of evidence in relation to proof and derivation has received rather less attention. I have suggested that it can be successfully accommodated within the framework of argumentation schemes. Specifically, appeals to evidence can be understood as C-schemes: argumentation schemes drawn from natural language reasoning that generally fall short of rigorous proof, but can sometimes be used rigorously. In this manner, mathematical arguments that rely on evidence to provide less than rigorous support for their conclusions can be understood as belonging to the same genus as mathematical proofs, but not the same species. Thereby the importance of rigour in proof is maintained, but without misleading the student into imagining that proof is somehow entirely alien from ordinary reasoning.

The framework of argumentation schemes, however, is itself a potentially valuable instrument for the mathematics educator.
Toulmin layouts, which are also drawn from the toolbox of argumentation theory, have lately been widely applied to mathematical reasoning by educational theorists
\cite[see, for example,][]{Knipping12}.
At least in their current form, argumentation schemes are a much more recent invention, and they have yet to receive much attention from mathematics educators: I am only aware of a few studies from one group of researchers specifically applying argumentation schemes to mathematics education \cite{Metaxas09,Metaxas15,Metaxas16,Koleza17}.
In addition, Aikaterini Konstantinidou and Fabrizio Macagno have discussed the application of argumentation schemes to science education \cite{Konstantinidou13,Macagno13c}.
They argue that argumentation schemes can have a particular value in `discovering the implicit beliefs affecting a student's learning process' \cite[235]{Macagno13c}. This is an important task in mathematics education too, and their study would seem to naturally generalize to mathematics education.

Perhaps the most exciting opportunity that argumentation schemes represent for mathematics education is that they are already widely implemented within what is sometimes called the `argument web'
\cite{Reed17}. In recent years, a growing number of software tools of increasing sophistication have been designed for the analysis and evaluation of argument.
These tools represent a potentially invaluable resource for the study of mathematical argument, but their application to mathematics has only just begun \cite{Pease17,Corneli18}. Again, I am unaware of any application of these resources in mathematics education: there is a growing body of work applying digital tools to mathematics education \cite[for example,][]{Modeste16,DurandGuerrier17}, but not the tools specific to argumentation. They represent an as yet untapped resource, of considerable scale and importance.

A final, more speculative benefit may accrue from an educational approach which emphasizes the continuity of mathematical argument and argument in other areas. Some authors have urged that mathematics educators promulgate the value of an education in mathematics as a source of the intellectual virtues and skills necessary for successful navigation of contemporary society \cite[recent examples include][]{Su17,Cheng18}. When mathematics is taught in a fashion that emphasizes its differences from everyday reasoning this can be a tough case to make. But argumentation schemes represent a plausible bridge between mathematical and other conversations; and mathematics may thereby provide an invaluable testbed for the acquisition and mastery of argumentation techniques of much wider application.

\section*{Acknowledgements}
I presented an earlier version of this paper at the interdisciplinary symposium on Mathematical Evidence and Argument held at the University of Bremen in 2017. I am grateful to the participants for their comments and particularly indebted to Christine Knipping and Eva M\"uller-Hill for their invitation and their hospitality in Bremen. I am also grateful to three anonymous referees for insightful and thorough comments.

\bibliographystyle{natbib-oup}\bibliography{ArgMath}\end{document}